\documentclass[12pt]{article}
\newtheorem{theorem}{Theorem}
\newtheorem{rem}{Remark}
\newtheorem{lemma}{Lemma}
\newtheorem{prop}{Proposition}
\newtheorem{cor}{Corollary}
\newtheorem{exmp}{Example}
\newtheorem{alg}{Algorithm}

\renewcommand{\epsilon}{\varepsilon}

\renewcommand{\hat}{\widehat}

\usepackage{epsfig}
\title{Energy Decay of Damped Systems}
\author{Kre\v simir Veseli\'c\thanks{Fernuniversit\" at
Hagen, Lehrgebiet
  Mathematische Physik, Postfach 940, D-58084 Hagen,
  Germany, e-mail:
    kresimir.veselic@fernuni-hagen.de.}}
\date{ }
\begin{document}
\maketitle
\begin{abstract} We present a new and simple bound for
the exponential decay of second order systems
using the spectral shift. This result is applied to
finite matrices as well as to partial differential
equations of Mathematical Physics. The type of the generated semigroup
is shown to be bounded by the upper real part of the numerical range of
the underlying quadratic operator pencil.
\end{abstract}
\section{Introduction and main estimate}\label{introduction}
In this note we consider the abstract second order system
\begin{equation}\label{system_MCK}
M\ddot{x} + C\dot{x} + Kx = 0
\end{equation}
where \(M,C,K\) are selfadjoint operators in
a Hilbert space \({\cal X}\) with the scalar product \(x^*y\)
linear in the second and antilinear in the first
variable.\footnote{All other conventions and notations will be taken
from \cite{Kato}.}
For simplicity, we assume that
the operators
\(M,K,C\) are positive definite and that
\(M,C\) are bounded.
The phase space formulation of (\ref{system_MCK}) reads
\begin{equation}\label{system_A}
\dot{y} = {\cal A}y,\quad {\cal A} =
\left[\begin{array}{rr}
0                  &  K^{1/2}M^{-1/2}           \\
- M^{-1/2}K^{1/2}  & -M^{-1/2}CM^{-1/2}         \\
\end{array}\right],
\end{equation}
\begin{equation}\label{system_A_1}
y =
\left[\begin{array}{r}
K^{1/2}x       \\
M^{1/2}\dot{x} \\
\end{array}\right],\quad
\end{equation}
with the solution
\[
y = e^{{\cal A}t}y^0.
\]
Thus, the square of the norm equals twice
the total energy of the system:
\[
\|y\|^2 = \|K^{1/2}x\|^2 + \dot{x}^*M\dot{x}.
\]
The operator \({\cal A}\) is readily seen to be maximal dissipative on
\[
{\cal D}({\cal A}) = {\cal D}(K^{1/2}) \oplus M^{1/2}{\cal D}(K^{1/2}).
\]
and thus the semigroup
\(e^{{\cal A}t}\) is contractive i.e. the energy of the system is a
non-increasing function of \(t\).\footnote{The reader primarily
interested in finite dimensional applications may skip
operator-theoretical details and just take \(M,C,K\) as
matrices. In this context 'positive' means 'positive
semidefinite'
and the expression  \(\|K^{1/2}x\|^2\) below may always
be read as
\(x^*Kx\).}
The resolvent of \({\cal A}\) is immediately seen to be
given by
\begin{equation} \label{resolvent}
 ({\cal A} - \lambda I)^{-1} =
 \left[ \matrix{-\frac{1}{\lambda}+\frac{1}
{\lambda}K^{1/2}K(\lambda)^{-1}K^{1/2} &  -K^{1/2}K(\lambda)^{-1}M^{1/2} \cr
& \cr
M^{1/2}K(\lambda)^{-1}K^{1/2} & -\lambda M^{1/2}K(\lambda)^{-1}M^{1/2} \cr }
\right]
\end{equation}
with
\begin{equation} \label{Kfamily}
 K(\lambda) = \lambda^2 M + \lambda C + K
\end{equation}
at least for those \(\lambda\) for which \(K(\lambda)\)
remains positive definite.
This formula is rigorous, if \(M\), \(C\), \(K\) are
all bounded, otherwise the terms
\(M^{1/2}K(\lambda)^{-1}K^{1/2}\) and
 \(K^{1/2}K(\lambda)^{-1}K^{1/2}\)
have to be replaced by their closures; the latter are
obviously everywhere defined and bounded.

Most existing works on the exponential decay estimate the
infimum of all
\(\beta\) for which
\begin{equation}\label{bound}
\|e^{At}\| \leq C_\beta e^{\beta t}.
\end{equation}

The infimum value \(\omega_0(A)\)  of all possible
\(\beta\) in (\ref{bound})
(the {\em type of the semigroup}
in the terminology of \cite{Kato}) is often equal to the maximal
real part
of the spectrum of \({\cal A}\).
In looking for this infimum usually little attention is
paid to the constant \(C_\beta\)
which may tend to infinity as
\(\beta\) approaches the infimum (see e.g. \cite{cox_zua},
\cite{batkai}, \cite{chen_zhou} ). Since this constant plays
a key role in controlling the finite-time behaviour of
the system
we are interested in a bound in which both  \(\beta\) and
\(C_\beta\)
are tried to be made
simply computable from the coefficients \(M\), \(C\), \(K\).
In fact, we obtain a family of estimates (\ref{bound}) for any
  \(\beta\) from the interval \((\gamma,0]\), where
 \begin{equation}\label{gamma00}
\gamma = \sup_{\stackrel{x \in {\cal D}(K^{1/2})}{x \neq 0}}
\Re\frac{- x^*Cx + \sqrt{(x^*Cx)^2-4x^*Mx\|K^{1/2}x\|^2}}{2x^*Mx},
\end{equation}
while \(C_\beta\) is an expression with similar ingredients.
The set
\[
W(K) = \{\lambda \in C; x^*K(\lambda)x = 0,\mbox{ for some unit } x\}
\]
is called {\em the numerical range} of the matrix pencil \(K(\lambda)\)
(cf. e.g \cite{markus_rodman}).
Our result implies, in fact,
\[
\omega_0({\cal A}) \leq  \gamma = \sup\Re W(K).
\]
The bound for \(\omega_0({\cal A})\) , obtained in \cite{batkai}
 reads in our notations
 \[
 \gamma_b = \max\left\{-\inf_x\frac{x^*Cx}{2x^*Mx},\:
 -\frac{1}{\sup_x\frac{x^*Cx}{x^*Kx} +
 2\sqrt{\sup_x\frac{x^*Kx}{x^*Mx}}}\right\}.
 \]
 The values \(\gamma\) and \(\gamma_b\) are not easy to compare
 in general. For any underdamped system i.e. whenever the expression
 under the square root in (\ref{gamma00}) is uniformly negative,
 we obviously have \(\gamma \leq \gamma_b\). For further comparisons
 see Sect. 2 below.

An estimate for \(C_\beta\)  was
obtained in \cite{chen_zhou}, Vol.~I, Ch.~6. for the wave
equation with distributed viscous damping. The bound obtained
there is much less explicit than ours; this actually made
impossible any comparison of the two.

As a by-product, we prove that
\(\omega_0({\cal A})\) is equal to the supremum of the real
part of the spectrum of \({\cal A}\) for ``partly overdamped systems''
i.e. for those for which \(2\gamma M + C\) is positive definite.
Another case in which this equality is shown to hold are the so-called
modally damped systems.

Our main tool will be the 'spectral shift substitution'
\[
x = e^{\mu t}z
\]
which gives rise to
a new phase space representation, equivalent to the previous one.
This yields a fairly
simple total-energy decay estimate. Applications are made
to both finite matrices and differential operators of
Mathematical
Physics.
The obtained estimate is shown to be (asymptotically)
attainable and
therefore in some weak sense optimal. On the other hand,
our estimate is void, if damping has a nontrivial nulspace.
While there are such systems with no exponential decay
at all, there are still many relevant cases whose decay is not
covered by our theory. In Sec.~2 we provide illustrating
examples.

The substitution \(x = e^{\mu t}z\) yields
\[
\dot{x} = \mu e^{\mu t}z + e^{\mu t}\dot{z}
\]
\[
\ddot{x} =
\mu^2 e^{\mu t}z + 2\mu e^{\mu t}\dot{z} + e^{\mu t}\ddot{z}
\]
and (\ref{system_MCK}) reads
\begin{equation}\label{Klnu}
M\ddot{z} + C(\mu)\dot{z} + K(\mu)z = 0
\end{equation}
with
\begin{equation}\label{Cmu}
C(\mu) = 2\mu M + C.
\end{equation}
As long as
\(C(\mu),K(\mu)\) stay positive definite this is equivalent to
the phase space representation
\begin{equation}\label{system_Amu}
\dot{w} = \hat{A}w,\quad \hat{A} =
\left[\begin{array}{rr}
0                  &  K(\mu)^{1/2}M^{-1/2}           \\
- M^{-1/2}K(\mu)^{1/2}  & -M^{-1/2}CM^{-1/2}         \\
\end{array}\right],
\end{equation}
\begin{equation}\label{system_Ahat_1}
w =
\left[\begin{array}{r}
K(\mu)^{1/2}z       \\
M^{1/2}\dot{z} \\
\end{array}\right],\quad
\end{equation}
with the solution
\[
w = e^{\hat{A}t}w^0.
\]
We now connect these two representations. We
have
\[
y_1 = K^{1/2}x = e^{\mu t}K^{1/2}z =
\mu e^{\mu t}K^{1/2}K(\mu)^{-1/2}w_1,
\]
\[
y_2 = M^{1/2}\dot{x} =
\mu e^{\mu t}M^{1/2}K(\mu)^{-1/2}w_1 +
e^{\mu t}w_2.
\]
Thus,
\begin{equation}\label{yLmu}
y =  e^{\mu t}{\cal L}(\mu)w,
\end{equation}
\begin{equation}\label{Lmu}
{\cal L}(\mu) =
\left[\begin{array}{rr}
K^{1/2}K(\mu)^{-1/2}        & 0 \\
\mu M^{1/2}K(\mu)^{-1/2}    & I \\
\end{array}\right].
\end{equation}
This, together with the evolution equations for \(y,w\)
gives
\begin{equation}\label{similarity}
 {\cal L}(\mu)\hat{A} = \left(A - \mu I\right){\cal L}(\mu).
\end{equation}
This yields the decay estimate
\begin{equation}\label{decay0}
\|e^{At}\| \leq \|{\cal L}(\mu)\| \|{\cal L}(\mu)^{-1}\|e^{\mu t}.
\end{equation}
The foregoing formal calculation is rigorous for
finite matrices.
Our general operator setting requires additional
justifications of these steps. Also, we need a more explicit
bound on the condition number which appears on the right hand
side of (\ref{decay0}). These two issues are the subject of
the following. We set
 \begin{equation}\label{p+-}
 p_\pm(x) =
 - \frac{x^*Cx}{2x^*Mx}\pm \sqrt{D(x)}
\end{equation}
with\footnote{By convention, the square root of a negative number
will have the positive imaginary part.}
\[
D(x) = \left(\frac{x^*Cx}{2x^*Mx}\right)^2 -\frac{\|K^{1/2}x\|^2}{2x^*Mx}
\]
and
\begin{equation}\label{gamma}
\gamma = \sup_{\stackrel{x \in {\cal D}(K^{1/2})}{x \neq 0}}
\Re p_+(x),
\end{equation}
\begin{prop} \label{infbound}
Let \(M,C,K\) be as defined above. Then (i)
\begin{equation}\label{gammalarger}
\gamma \geq - \inf_x \frac{x^*Cx}{2x^*Mx},
\end{equation}
moreover,  \(\gamma\) is the infimum of all \(\mu\) for which
both \(K(\mu)\) and \(2\mu M + C\) are positive definite.
\end{prop}
{\bf Proof.} The relation (\ref{gammalarger}) is obvious. To prove the
second assertion note that for any real \(\mu\) the selfadjoint operator
\(K(\mu)\) is generated by the form
\[
h_\mu(x,y) = \mu^2x^*My + \mu x^*Cy + (K^{1/2}x)^*K^{1/2}y
\]
which is symmetric, bounded from below and closed on \({\cal D}(K^{1/2})\).
Take any \(\mu > \gamma\). Then
\[
h_\mu(x,x) = x^*x(\mu - p_-(x))(\mu - p_+(x))
 \geq x^*x(\mu - \Re p_+(x))^2
\geq x^*x(\mu - \gamma)^2
\]
(note that \( p_-(x)\leq p_+(x)\) whenever \(D(x)  \geq 0\)).
 Thus,
\(K(\mu)\) is positive definite. By (\ref{gammalarger}) \(2\mu M + C\)
is positive definite also.

Conversely, suppose that \(2\gamma M + C\) is positive definite;
there is a
sequence \(x_n\) of unit vectors such that
\[
\Re p_+(x_n) = - \frac{x_n^*Cx_n}{2x_n^*Mx_n}+ \Re \sqrt{D(x_n)}
\rightarrow \gamma, \quad n \rightarrow \infty.
\]
Since \(2\gamma M + C\) is positive definite it follows
\[
\lim\inf_{n \rightarrow \infty} \Re\sqrt{D(x_n)}  > 0.
\]
Thus, for \(n\) sufficiently large \(D(x_n) > 0\), hence
\[
p_+(x_n) = \Re p_+(x_n) \rightarrow \gamma,
\quad n \rightarrow \infty
\]
and (\(K(\gamma)\) is positive semidefinite)
\[
h_\mu(x_n,x_n) = \|K(\gamma)^{1/2}x_n\|^2 =
(\gamma - p_-(x_n))(\gamma - p_+(x_n))
\rightarrow 0
\]
since \(p_-(x_n)\) is obviously bounded. Thus,  \(K(\gamma)\)
has not a bounded inverse, and is not
positive definite. Q.E.D.

\begin{theorem}\label{theorem}
Let \(M,C,K\) be as defined above. Assume that
\(\gamma < \mu \leq 0\) holds.
Then (i) the operator \({\cal L}(\mu)\)
together with its inverse
is bounded and everywhere defined. (ii) Both operators
leave the subspace \({\cal D}(A)\) invariant and
(\ref{similarity}) holds. (iii) The estimate  (\ref{bound})
holds with
\begin{equation}\label{final_est}
\beta = \mu, \quad
C_\beta =  \|{\cal L}(\mu)\| \|{\cal L}(\mu)^{-1}\|.
\end{equation}
(iv) The type of the semigroup is bounded by
 \begin{equation}\label{type_bound}
\omega_0({\cal A}) \leq \gamma.
\end{equation}
\end{theorem}
{\bf Proof.}
Since \(\gamma\) is negative then for \(-\gamma < \mu \leq 0\)
both \(K(\mu)\) and \(C(\mu)\) remain positive definite.
The boundedness of \(M,C\) implies (cf. \cite{Kato})
\[
{\cal D}(K(\mu)^{1/2})= {\cal D}(K^{1/2}),\quad \mbox{and}\quad
{\cal D}(K(\mu))= {\cal D}(K),\quad
0 >\mu > \gamma,
\]
so both (\ref{Lmu}) and its inverse
\begin{equation}\label{Lmu-1}
{\cal L}(\mu)^{-1} =
\left[\begin{array}{rr}
K(\mu)^{1/2}K^{-1/2}          & 0 \\
-\mu M^{1/2}K^{-1/2}          & I \\
\end{array}\right]
\end{equation}
 are everywhere
defined and bounded (a straightforward calculation
shows that the one is the inverse of the other).
This also shows that \({\cal D}(A) = {\cal D}(\hat{A})\).

Take a vector from \({\cal D}(A)\), that is, \(x \in {\cal D}(K^{1/2})\),
\(y \in M^{1/2}{\cal D}(K^{1/2})\) and set
\[
{\cal L}(\mu)
\left[\begin{array}{c}
x  \\
y  \\
\end{array}\right]
=
\left[\begin{array}{c}
u  \\
v  \\
\end{array}\right].
\]
Here \(x = K^{1/2}K(\mu)^{-1/2}x\) is from \({\cal D}(K^{1/2})\)
because \(K(\mu)^{-1/2}x \in {\cal D}(K(\mu) =  {\cal D}(K)\).
Also \(v = \mu M^{1/2}K(\mu)^{-1/2}x + y\) is from
\(M^{1/2} {\cal D}(K^{1/2})\) because \(y\) is such and
\(K(\mu)^{-1/2}x \in  {\cal D}(K)\). Thus, \({\cal L}(\mu)\)
leaves   \({\cal D}(A)\) invariant. The proof for  \({\cal L}(\mu)^{-1}\)
is similar. Now the relation (\ref{similarity}) can be directly
verified on any vector from  \({\cal D}(A)\). The estimate
(\ref{final_est}) is now obvious.
Q.E.D.\\

\begin{cor} \label{abscissa} If
\[
\gamma > - \inf_x \frac{x^*Cx}{2x^*Mx},
\]
then
\[
\gamma = \omega_0({\cal A}) = \sup\sigma({\cal A}).
\]
\end{cor}
 {\bf Proof.} Under our assumption, as in the proof of
 Proposition \ref{infbound}, it follows that \(K(\gamma)\)
 is not boundedly invertible while \(K(\lambda)\) is positive
 definite for any \(\lambda > \gamma\). Now, the  \(2,2\)-block
 in the resolvent matrix in  (\ref{resolvent}) gets unbounded
 for  \(\lambda = \gamma\). Thus, all \(\lambda\) with
 \(\Re\lambda > \gamma\) belong to the resolvent set of
  \({\cal A}\) whereas \(\gamma\) belongs to its spectrum,
  hence \(\gamma = \sup\Re\sigma({\cal A})\).
  From Theorem \ref{theorem} it follows
   \(\omega_0({\cal A})\leq \gamma\)  while
   \(\omega_0({\cal A}) \geq \sup\sigma({\cal A})\)
  is a general fact. So, the assertion follows. Q.E.D.

The systems covered by the corollary above may be called
'partially overdamped' in  a sense similar to that
introduced in \cite{barston}.

We next give some bounds for the condition number appearing in
(\ref{final_est}). Note that both \({\cal L(\mu)}\) and
\({\cal L(\mu)}^{-1}\) are of the type
\begin{equation}\label{type}
{\cal L} =
\left[\begin{array}{cc}
A   &   0   \\
B   &   I   \\
\end{array}\right].
\end{equation}

\begin{lemma} \label{LemmaL}
Let \(A\),\(B\) be bounded operators and
Then \(\|{\cal L}\|\) from (\ref{type}) is bounded by any
of the numbers
\begin{equation}\label{quantities}
\sqrt{1 + \|A^*A + B^*B\|},\quad \max\{\|A\|,1 + \|B\|\}
,\quad
 \max\{\|A\| + \|B\|,1\}.
\end{equation}
\end{lemma}
{\bf Proof.} We have
\[
{\cal L}^*{\cal L} =
\left[\begin{array}{cc}
A^*A + B^*B    &   B^* \\
B              &   I   \\
\end{array}\right], \quad
{\cal L}{\cal L}^* =
\left[\begin{array}{cc}
AA^*     &   AB^*      \\
BA^*     & BB^* +  I   \\
\end{array}\right].
\]
If \(B = 0\) the assertion is trivial. If \(B \neq 0\) the number
\(a = \|{\cal L}\|^2 =  \|{\cal L}{\cal L}^*\| =  \|{\cal L}^*{\cal L}\|
> 1\) belongs to the spectrum of \({\cal L}^*{\cal L}\) and there
exist sequances \(x_n\), \(y_n\) with \(\|x_n\|^2 + \|y_n\|^2 = 1\)
and
\begin{eqnarray}\label{block1}
(A^*A + B^*B)x_n + B^*y_n - ax_n &  \rightarrow & 0 \\ \label{block2}
Bx_n + y_n - ay_n                &  \rightarrow & 0 \\
\end{eqnarray}
for \(n \rightarrow 0\). Hence
\begin{equation}\label{quadratic}
\left(A^*A + B^*B + \frac{B^*B}{a - 1} - a\right)x_n  \rightarrow 0
\end{equation}
By (\ref{block2}) the sequence \(x_n\) does not converge to \(0\)
and without loss of generality we may assume that in (\ref{block1})
\(\|x_n\| = 1\) holds. Then
\[
a^2 - a((A^*A + B^*B + I)x_n,x_n) + (A^*Ax_n,x_n)  \rightarrow 0.
\]
Hence
\[
a \leq \sup_{\|x\| = 1}
\frac{((A^*A + B^*B + I)x,x) + \sqrt{(A^*A + B^*B + I)x,x)^2 - 4(x,x)(A^*Ax,x)}}
{2}
\]
\[
\leq
 \sup_{\|x\| = 1}((A^*A + B^*B + I)x,x).
\]
The other two bounds in (\ref{quantities}) are obvious. Q.E.D.

\begin{rem}{\em
(i)  Ours is, in fact a family of bounds depending
on the parameter \(\mu\) from the interval \((\gamma,0]\).
An optimal bound would be obtained as the infimum over all
of them.

(ii) By continuity, the bound (\ref{final_est})
remains valid even for \(\mu = \gamma\), if
\(K(\gamma)\) is positive definite. In this case we have
\(\omega_0({\cal A}) < \gamma\) i.e. or bound is not optimal.
Indeed, for some \(\mu < \gamma\) the operator \(K(\mu)\) will
still be positive definite and  \(\hat{A}\) will still generate
a uniformly bounded semigroup (\cite{Kato}, Ch.~IX, Th.~2.1)
and (\ref{similarity}), (\ref{decay0}) will still be valid.

(iii) Our main estimate does not contain the norms
of the operators \(M,C\) or their inverses
and indeed it would certainly
hold in much more singular cases, but then some additional
regularity conditions would be needed e.g. the boundedness
of the operators \(C^{1/2}K^{-1/2}\),
\(C^{-1/2}M^{1/2}\) and the like thus requiring more technical
proofs. The same is valid, if we would admit more
general damping operator \(C\) by merely asking it to be accretive.

(iv) The expression (\ref{gamma}) for \(\gamma\) is neat but not
easily computable, even in the case of finite matrices. As suggested
by Proposition  \ref{infbound} a simple viable
method to determine \(\gamma\) numerically would run as follows
\begin{itemize}
\item Find \(\gamma_0 = -\inf \frac{x^*Cx}{2x^*Mx}\) by finding the lowest
eigenvalue of the matrix pencil \(2\lambda M - C\).
\item If \(\gamma_0 \approx 0\) halt, no bound available.
\item If \(\gamma_0 < 0\) and \(K(\gamma_0)\) is positive definite
then \(\gamma =\gamma_0\).
\item If \(\gamma_0 < 0\) and \(K(\gamma_0)\) is not positive definite
then seek \(\gamma\) by bisection in the interval \((\gamma_0,0]\).

\end{itemize}
}
\end{rem}
 \section{Examples}\label{examples}
Our first example is the one-dimensional system with
the two-dimensional phase space matrix
\begin{equation}\label{A2x2}
A =
\left[\begin{array}{rr}
0        &   k \\
-k       &  -d  \\
\end{array}\right],
\quad k > 0,\quad d > 0.
\end{equation}
A straightforward, if a bit tedious, computation gives
\begin{equation}\label{expA2x2}
e^{At} =
\left[\begin{array}{cc}
\cos\delta t + d(\sin\delta t)/(2\delta) & k(\sin\delta t)/\delta \\
-k(\sin\delta t)/\delta & \cos\delta t - d(\sin\delta t)/(2\delta)\\
\end{array}\right]e^{-dt/2}.
\end{equation}
with
\[
\delta =  \sqrt{4k^2 - d^2}/2.
\]
The above formula is valid for \(\delta\)
both positive and negative (in the latter case there
is a real expression by means of hyperbolic functions)
whereas the formula for \(4k^2 - d^2 = 0\) is obtained
taking the limit \(k \rightarrow d/2\), thus reading
\begin{equation}\label{expA2x20}
e^{At} =
\left[\begin{array}{cc}
1 + dt/2  & dt/2 \\
 -dt/2 & 1 -dt/2\\
\end{array}\right]e^{-dt/2}.
\end{equation}
Here
\[
{\cal L}(\mu) =
\left[\begin{array}{cc}
\frac{k}{\sqrt{\mu^2 + \mu d + k^2}}   &   0 \\
\frac{\mu}{\sqrt{\mu^2 + \mu d + k^2}} &   1 \\
\end{array}\right]
\]
and our Theorem above yields
\begin{equation}\label{omega0x20}
\gamma = \Re\frac{-d + \sqrt{d^2 - 4k^2}}{2}
\end{equation}
Here the right hand side is the largest real part of the
spectrum and therefore the inequality (\ref{type_bound}) is,
in fact,
an equality i.e. our estimate for
\(\omega_0({\cal A})\) is sharp
for all possible values of \(k,d\). The bound
from \cite{batkai} reads (in our notations)
\[
\gamma_b =
\max\left\{-\frac{k^2}{d + 2k},-\frac{d}{2}\right\}.
\]
As a straightforward calculation shows we have here
\(\gamma_b = \gamma\) for \(d \leq (\sqrt{3} - 1)k\) and
\(\gamma_b > \gamma\) for \(d > (\sqrt{3} - 1)k\).

\begin{figure}[ht]
\begin{center}
\epsfig{file=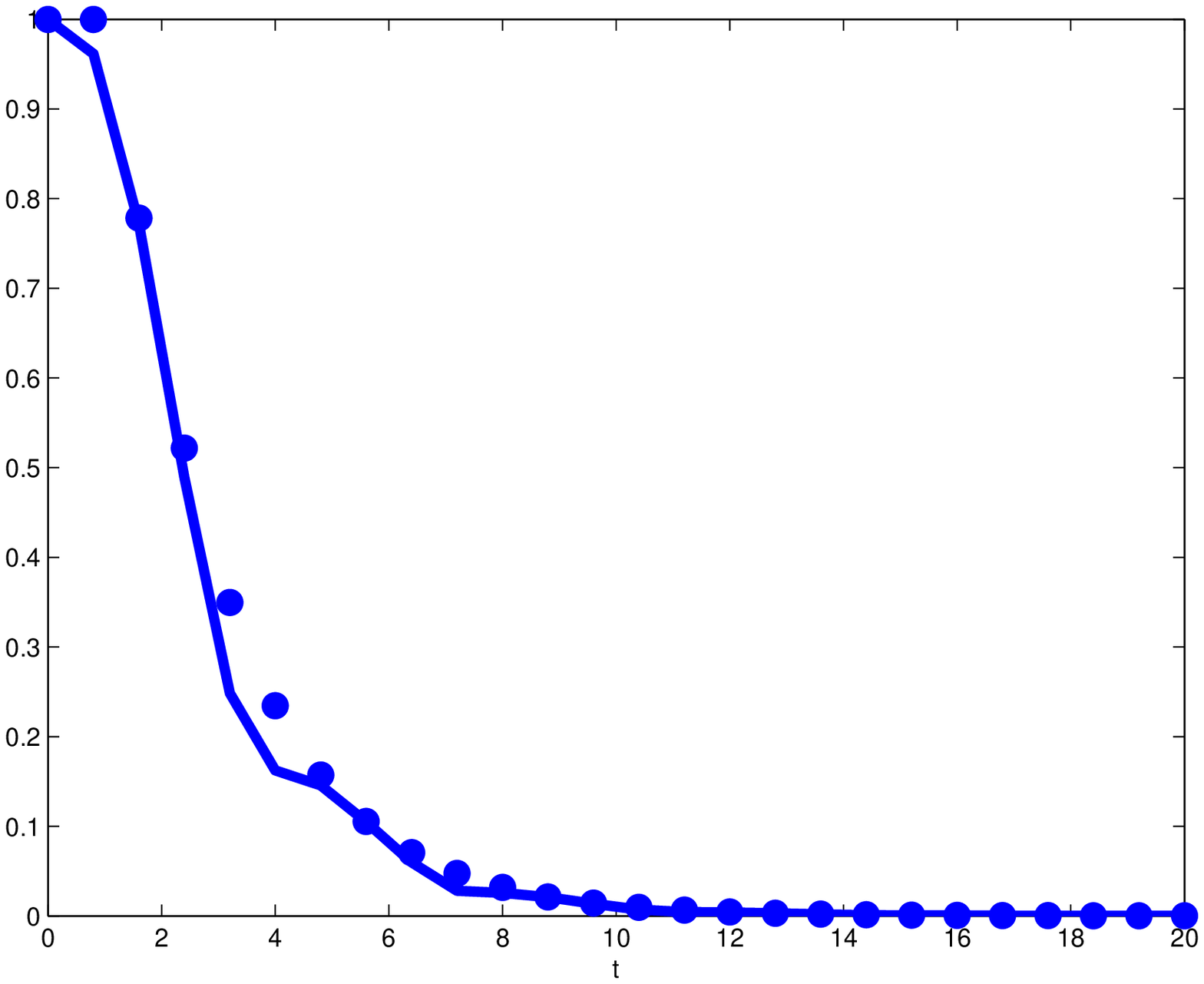,scale=0.6}
\end{center}
\caption{\label{kd1} \(k/d = 1\)}
\end{figure}\hfill \\

\begin{figure}[ht]
\begin{center}
\epsfig{file=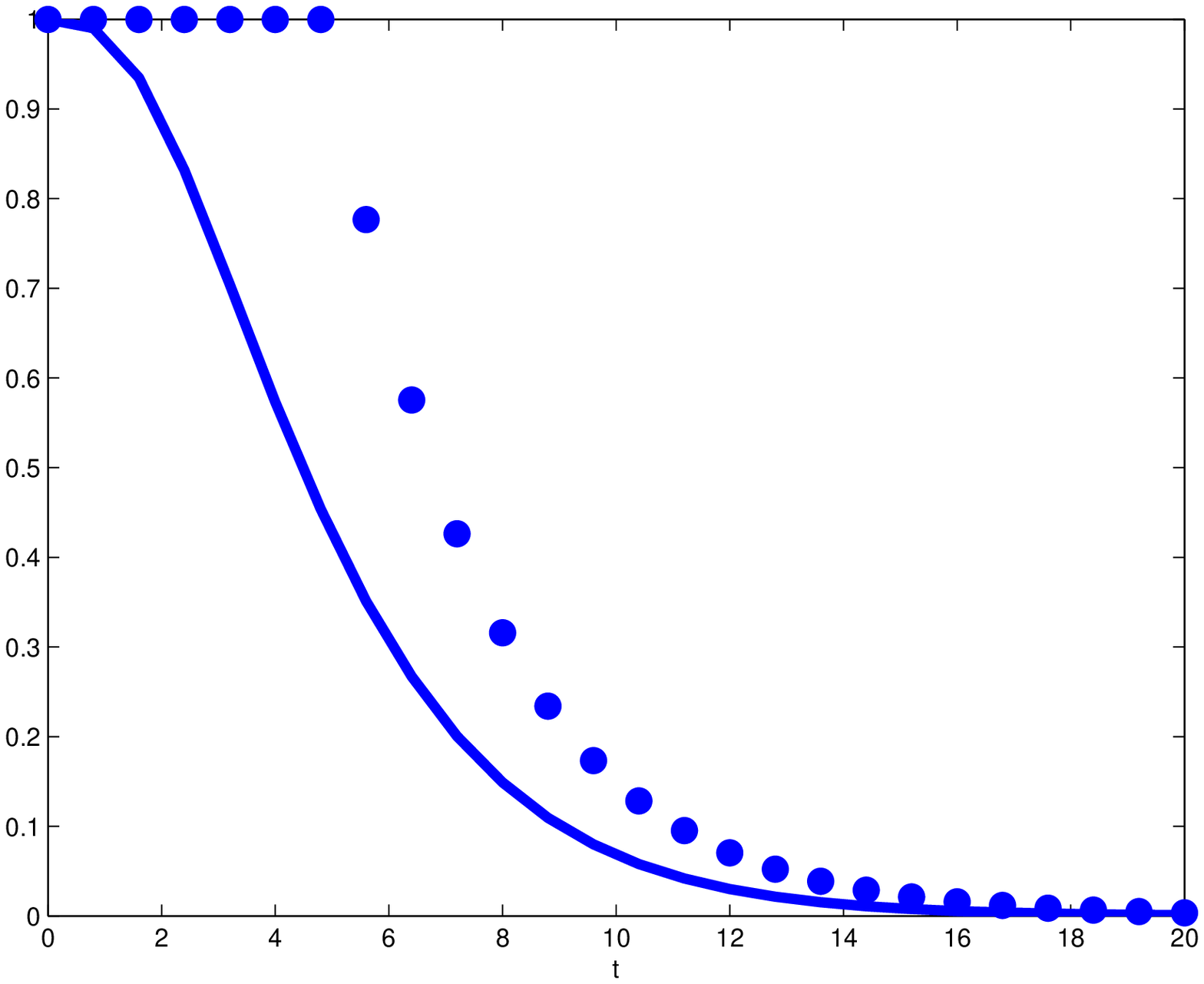,scale=0.6}
\end{center}
\caption{\label{kd2} \(k/d = 0.5\)}
\end{figure}\hfill \\

On Figs. \ref{kd1}, \ref{kd2} we display our bound
(dotted line) with the true norm (solid line) as
functions of \(t\) (our bound is obtained by taking
the minimum over four equidistant values of \(\mu\)).
We see that our constant \(C_\beta\) gets
pessimistic  for small values of \(k/d\). On the other
hand, for \(k/d \rightarrow \infty\) we may
take \(\mu = d/2\),
\({\cal L}(\mu) \rightarrow I\) and
(\ref{final_est}) reads asymptotically
\[
\|e^{At}\| \leq e^{dt/2}.
\]
whereas the right hand of  (\ref{expA2x2}) behaves as
\[
\left[\begin{array}{cc}
\cos kt  &   \sin kt \\
-\sin kt &   \cos kt \\
\end{array}\right]e^{-dt/2}.
\]
Thus, or bound (\ref{final_est}) is asymptotically sharp.

Any modally damped system i.e. a system in which
\(M,K,C\) satisfy the relation
\[
CM^{-1}K = KM^{-1}M
\]
is easily seen to be unitarily equivalent to an orthogonal
sum of matrices of type (\ref{A2x2}), a bound is
obtained as the maximum over all of them.  The type
\(\omega_0({\cal A})\) is the maximum of the values in
(\ref{omega0x20}).\footnote{This remark is literally true
if the system has discrete spectrum (e.g., if \(K\) has
a compact inverse). Otherwise we would have direct
integrals.} In particular, \(\omega_0({\cal A})\) is equal to
\(\Re\sigma({\cal A})\).

As a second example take the wave equation in a bounded
domain \(\Omega \subset R^n\)
\begin{equation}\label{waveq}
w_{tt}(x,t) + c(x)w_t(x,t) +\Delta w(x,t) = 0,
\end{equation}
with the boundary condition
\[
w(x,t) = 0 \quad x \in \partial\Omega.
\]
The function \(c\) is supposed to have finite positive
minimum and maximum.

We will estimate e.g. the second of the bounds in
(\ref{quantities}). Taking first
\[
A = K^{1/2}K(\mu)^{-1/2},\quad
B = \mu M^{1/2}K(\mu)^{-1/2}
\]
 in \({\cal X} = L_2(R)\) we have
\[
\|A\|^2 =
\frac{1}{1 + \inf_u\frac{u^*(\mu^2 + \mu c(\cdot))u}
{\|\nabla u\|^2}} =
\frac{1}{1 + \lambda_1},
\]
where \(\lambda_1\) is the lowest eigenvalue \(\lambda\)
of the boundary value problem
\[
(\mu^2 + \mu c(x))u = -\lambda\Delta u,\quad
u|_{\partial \Omega} = 0.
\]
Since under our assumptions \(\mu^2 + \mu c(x)\)
is negative definite we have
\[
\lambda_1 < 0,\quad 1 + \lambda_1 > 0.
\]
Furthermore
\[
\|B\|^2 =
\sup_u\frac{\mu^2u^*u}
{u^*(\mu^2 + \mu c(\cdot) - \Delta)u} =
\frac{1}{\mu^2 + \inf_u\frac{u^*(\mu c(\cdot) - \Delta)u}
{u^*u}}
\]
\[
= \frac{1}{\mu^2 + \lambda_2},
\]
where \(\lambda_2\) is the lowest eigenvalue \(\lambda\)
of the boundary value problem
\[
(\mu c(x) - \Delta)u = \lambda u,\quad
u|_{\partial \Omega} = 0.
\]
Taking next
\[
A = K(\mu)^{1/2}K^{-1/2},\quad
B = -\mu M^{1/2}K^{-1/2},
\]
we have \(\|A\| \leq 1\) and
\[
\|B\|^2 =
\frac{1}
{\lambda_3}
\]
 where
\(\lambda_3\) is the lowest eigenvalue \(\lambda\)
of the boundary value problem
\[
 \Delta u = -\lambda \mu u,\quad
u|_{\partial \Omega} = 0.
\]
Altogether
\[
C_\mu = \|{\cal L}(\mu)\|\|{\cal L}^{-1}(\mu)\|
\leq
\max\left\{\frac{1}{\sqrt{1 + \lambda_1}}, 1 +
\frac{1}{\sqrt{\mu^2 + \lambda_2}}\right\}
\left(1 + \frac{1}{\sqrt{\lambda_3}}\right),
\]
where \(\lambda_1, \lambda_2, \lambda_3\) are obtained above.
Thus, our bound for \(C_\mu\) is obtained from the extremal
eigenvalues of some selfadjoint elliptic boundary value problems
involving \(M,C,K\).
If \(c(x)\) is constant then all these boundary value problems
reduce to \(-\Delta u = \lambda u\). This system is also
modally damped, so the corresponding estimates are applicable
here, too.

\end{document}